\documentstyle[12pt]{article}

\newcounter{const}

\newcommand{\argm}{\mathop{\rm argmax}\limits}

\begin{document}



\begin{center}
\vspace{2mm}
{\bf Exact Constants in the Rosenthal Moment Inequalities for Sums of 
independent centered Random Variables.}\\
\vspace{4mm}

{\bf Naimark B.} \footnote{Department of Mathematics and Computer science, Open 
University of Israel, Raanana, 43107, Israel; 108 Ravutski street, POB 808.
 E - mail: n-boris@netvision.co.il} 
{\bf Ostrovsky E.} \footnote{Corresponding author. Department of Mathematics,
 Ben - Gurion University, Beer - Sheva,  Israel. E - mail: galaostr@cs.bgu.ac.il}
\end{center}

\vspace{3mm}

{\bf Abstract.}
We study the exact constants in the moment inequalities for sums of centered 
independent random variables: improve their asymptotics, low and upper bounds,
calculate more exact asymptotics, elaborate the numerical algorithm for their 
calculation, study the class of smoothing etc.


\vspace{3mm}







\vspace{1mm}
\vspace{1mm}





\vspace{2mm}


\vspace{2mm}

\section {Introduction. Statement of problem.} 

 Let $ \ p = const \ge 2, \  \{ \xi(i) \}, \ i = 1,2,\ldots, n $ be a sequence of independent centered: 
$ {\bf E} \xi(i) = 0 $ random variables belonging to the space $ L_p, $ i.e. such that 
$$
\forall i \  \  ||\xi(i)||_p \stackrel{def}{=}  {\bf E}^{1/p} |\xi(i)|^p < \infty. \eqno(0)
$$
 We denote $ \sum a(i) = \sum_{i=1}^n a(i), \ L(p) = C^p(p), $ where 

$$
C(p) = \sup_{ \{\xi(i) \} } \sup_{n} \frac{||\sum \xi(i)||_p  }
{\max \left(||\sum \xi(i)||_2, \left(\sum ||\xi(i)||^p_p \right)^{1/p} \right) }, \eqno(1)
$$
( $ "C" $ denotes the centered case), 
 where the external $ "\sup" $ is calculated over all the sequences of independent centered
 random variables satisfies the condition (0). \par
 In the case if in (1) all the variables $ \{ \xi(i) \} $ are symmetrically distributed,
independent, $ ||\xi(i)||_p < \infty, $ we will denote the correspondent 
constants (more exactly, functions of $ p) $ 
 $ S(p) \ ( "S" (\cdot) $ denotes the symmetrical case) 
instead $ C(p) $ and  $ \ K(p)\stackrel{def}{=} S^p(p) $ instead $ L(p). $  It 
is obvious that $ S(p) \le C(p), \ K(p) \le L(p). $ In the article [12] is proved that
$ C(p) \le 2 S(p), \ L(p) \le 2^p \ K(p).$  \par
 The constant $ C(p), S(p) $ are called the exact constants in the moment inequalities for the 
sums of independent random variables  and 
 play very important role in the classical theory of probability ([1], 522 - 523,
[2], p. 63;) theory of probability on the Banach spaces [4], in 
the statistics and theory of Monte - Carlo method ([19], section 5) etc. \par
 There are many publications on the behavior of constants $ C(p), S(p) $ at $ \ p \to \infty. $
The first estimations are obtained in [6]; Rosenthal [7] proved in fact that  
$ C(p) \le C^p_1;  \ C_1 = const > 1; $
 here and further $ C_j, \ j = 1,2,\ldots $ are some positive finite
{\it absolute} constants, $ \log = \ln. $ In the article [17] is proved that 
$ C(p) \le 9.6 \ p / \log p. $ In the works [8], [9], [22] are obtained the 
{\it non - asymptotical } bide - sides estimations  for $ S(p): $

$$
(e \sqrt{2})^{-1} \ p /\log p \le S(p) \le 7.35 \ p /\log p, \ p \ge 2, \eqno(2)
$$
and there are some moment estimations for the sums independent nonnegative random variables.
 See also Latala [12], Utev [17], [18]; Pinelis and Utev
 [20] and so one.\par
 In the articles of Ibragimov R. and Sharachmedov Sh. [10], [11] and Utev  [17], [18]
is obtained the {\it explicit } formula for  $ S(p): \ S(2) = 1; $ at $ p \in (2,4] $ 

$$
S(p) = \left( 1 + \sqrt{ \frac{2^p}{\pi} } \Gamma \left( \frac{p+1}{2} \right)  \right)^{1/p};
$$

$$
p \ge 4 \ \Rightarrow S(p) = ||\tau_1 - \tau_2||_p, \eqno(3)
$$
where $ \Gamma(\cdot) $ is the Gamma function, a random variables $ \tau_j $ are independent and 
have the Poisson distribution with parameters 0.5: $ \ {\bf E} \tau_j = {\bf D} \tau_j = 1/2. $ \par
 As a consequence was obtained that at $ p \to \infty $ 

$$
S(p) = \frac{p}{e \cdot \log p} \left(1 + o \left(\frac{\log^2 \log p}{\log p} \right)  \right).
$$

 In the article [11] is obtained the following representation for the values $ L(2m), m = 2, 4, 6,
\ldots: $
$$
C^{2m}(2m) = L(2m) = {\bf E}(\theta - 1)^{2m} = e^{-1} \sum_{n=0}^{\infty} (n-1)^{2m} /n!,
$$
where the random variable $ \theta $ has the Poisson distribution with parameter 1, and there 
is a hypothesis that for all the values $ p \ge 4 \ C^p(p) = L(p). $ \par
 We will denote also for all values $ p \ge 4 $

$$
L(p) = {\bf E}|\theta -1|^p = e^{-1} \sum_{n=0}^{\infty} |n-1|^p/n!, \ \ G(p) = L^{1/p}(p).
\eqno(4)
$$

 In the report [13] are obtained the estimations for $ C(p) $ in the case if the 
sequence $ \{\xi(i) \} $ is the sequence 
of martingale - differences, in the article [18] there are some generalizations 
for weakly dependent random variables $ \{\xi(i) \}. $ \par
 {\bf In this article we improve the bide - side estimations and asymptotics for $ S(p), 
\ G(p) $ at $ p \to \infty, $
find the exact boundaries for the different approximation of $ S(p), \ G(p); $
describe the algorithm for the numerical calculation of $ K(p), \ L(p); $ 
study the analytical properties  of $ K(p), \ L(p) $ etc.}\par 
 Notice that there are many other statements of this problem: for the nonnegative variables
[8], [12]; for the Hilbert space valued variables [18] etc. \par 

\section{Main results.}

 Let us introduce the following functions at $ p \ge 4: \ g(p) = p/(e \log p), \ 
\delta(p) = 1/\log p, \ \Delta(p) = \log \log p/\log p; $ 

$$
  h(p) = g(p) \left(1+ \Delta(p) + \Delta^2(p) \right) =  
$$

$$
 [p/(e \ \log p)] \cdot  \left(1+ \log \log p/\log p + (\log \log p/\log p)^2 \right);
$$

$$
I_n(z) = 2^{-n} \sum_{k=0}^{\infty} 4^{-k} z^k /(k! \ (n+k)!)
$$
is the usually modified Bessel's function of order $ n; $ 

$$
W(p) = (2/e) \sum_{n=1}^{\infty} n^p I_n(1);
$$

$$
B(p) = e^{-1} \sum_{n=1}^{\infty} n^p/n!, \ p >0; \ B(0) = 1;
$$
$ B(p) $ are the well - known Bell's numbers; $ B(p) = {\bf E} \tau^p, $ where the random 
variable $ \tau $ has the usually Poisson distribution with parameter 1: $ {\bf E} \tau =  
{\bf D} \tau = 1. $ \par
 The generalized Bell's function $ B(a,p;z) $ may be defined as

$$
B(a,p,z) = \sum_{n=0}^{\infty} \frac{|n-a|^p \ z^n }{e \cdot n!}.
$$
 For example, $ B(0,p,1) = B(p). $ \par

{\bf Theorem 1.}

$$
1 = \inf_{p \ge 4} G(p)/g(p) < \sup_{p \ge 4} G(p)/g(p) = C_3,  \eqno(5a)
$$

{\it where}
$$
C_3 = \sup_{p \ge 4} B^{1/p}(1,p,1)/g(p) = G(C_4)/g(C_4) \approx 1.77638,
$$

$$
C_4 = \argm_{p \ge 4} B^{1/p}(1,p,1)/g(p) \approx 33.4610; 
$$
{\it (The equality} $ C_3 \approx 1.77638 $ {\it means that } $ |C_3 - 1.77638| 
\le 5 \cdot 10^{-6}); $ 

$$
 1 = \inf_{p =4,6,8 \ldots} C(p) /g(p) < \sup_{p=4,6,8 \ldots} C(p)/g(p) = C_5, \eqno(5b)
$$

{\it where } $ C_5 = $

$$
 \inf_{p=4,6,8,\ldots} B^{1/p}(1,p,1)/g(p) = G(C_6)/g(C_6) \approx 1.77637, \ C_6 = 34;
$$

$$
1 = \inf_{p \ge 15} G(p)/h(p) < \sup_{p \ge 15}G(p)/h(p) = G(C_8)/h(C_8) = C_{7}, \eqno(5c)
$$

{\it where}

$$
C_7 = \sup_{p \ge 15} B^{1/p}(1,p,1)/h(p) \approx 1.2054,
$$

$$
C_8 = \argm_{p \ge 15} B^{1/p}(1,p,1)/h(p) \approx 71.430;
$$

$$
1= \inf_{p \ge 4} S(p)/g(p) < \sup_{p \ge 4} S(p)/g(p) = C_9, \eqno(6a) 
$$

{\it where }

$$
C_9 = \sup_{p \ge 4} W^{1/p}(p)/g(p) = S(C_{10})/g(C_{10}) \approx 1.53572, 
$$

$$
 C_{10} = \argm_{p \ge 4} W^{1/p}(p)/g(p)  \approx 22.311; 
$$

$$
1 = \inf_{p \ge 15} S(p)/h(p) < \sup_{p \ge 15} S(p)/h(p) = S(C_{12})/h(C_{12})= C_{11}, \eqno(6b)
$$

{\it where}

$$
C_{11}= \sup_{p \ge 15} W^{1/p}(p)/h(p) \approx 1.03734,
$$

$$
C_{12} = \argm_{p \ge 15} W^{1/p}(p)/h(p) \approx 138.149;
$$

$$
1 = \inf_{p = 16,18,20, \ldots} C(p)/h(p) < \sup_{p=16,18.20,\ldots} C(p)/h(p) = 
$$

$$
C(72)/h(72) = \sup_{p = 16,18,20,\ldots} B^{1/p}(1,p,1)/h(p) \approx 1.2053. \eqno(6c)
$$

(We choose the value 15 as long as the function $ \log \log p/\log p $ monotonically 
decreases for the values $ p \ge \exp(e) \approx 15.15426 ). $ \par
 Notice than our estimations and constants 
 (5a, 5b, 5c) and (6a, 6b,6c) are exact and improve the constants and estimations of 
Rosenthal [7]; Johnson, Schechtman, Zinn [8]; Ibragimov, Sharachmedov [10],[11]; Latala [12];
 Utev [17], [18] etc. 
For example, $ 1/(1/\sqrt{2}) \approx 1.41421, \ 7.35 e/C_3 
\approx 11.2472. $ \par

\vspace{2mm}

{\bf Theorem 2. }  {\it  At } $ \ p \to \infty  \ G(p) = [p/(e \cdot \log p)] \times $  
$$
 \left(1 + \frac{ \log \log p}{\log p} + \frac{1}{\log p} + 
\frac{\log^2 \log p}{\log^2 p} + \frac{ \log \log p}{\log^2 p} (1 + o(1)) \right); \eqno(7a)
$$
 $ S(p) = [p/(e \cdot \log p)] \times $

$$
 \left( 1 + \frac{\log \log p}{\log p} + \frac{1 - \log 2}{\log p} + \frac{\log^2 \log p}
{\log^2 p} + o \left( \frac{\log \log p}{\log^2 p }  \right) \right). \eqno(7b)
$$

 Let us denote for the values $ p \ge 4 $
 by $ N = N(p), \ M = M(p) $  the (unique) solutions of
equations

$$
M(p) \log M(p) = p, \ \ N(p) \log(2 N(p)) = p, \eqno(8)
$$
such that $ N(p) = 0.5 M(2p). $ \par

{\bf Theorem 3. } {\it At } $ p \to \infty, \ m = 2,3,4, \ldots \to \infty $

$$
G(p) = M(p)^{1-M(p)/p} \ \exp(M(p)/p) \ (1+O(\log p/p)), \eqno(9a)
$$

$$
C(2m) = M^{1-M(2m)/2m}(2m) \ \exp(M(2m)/(2m)) \ (1+O(\log m/m)),
$$

$$
S(p) = N \ (e/2N)^{N/p} \ \ (1 + O(\log p/p)). \eqno(9b). 
$$

 Denote by $ s(n,r) $ the usually Stirling's numbers of a second kind appeared in the 
combinatorics ([14], p. 117):

$$
x^n = \sum_{r=0}^n s(n,r) x_{(r)}; \ x_{(r)} \stackrel{def}{=} x(x-1)(x-2)\ldots (x-r+1), \ 
x_{(0)} = 1.
$$

{\bf Theorem 4.} {\it Let } $ p $ {\it be even: } $ p = 2m, \ m=2,3,4, \ldots. $ {\it Then }

$$
K(2m)=\sum_{l=0}^{2m} (-1)^l {2m \choose l} \sum_{q=0}^{2m-l} \sum_{r=0}^l 2^{-r-q} s(2m-l,q)
s(l,r), \eqno(10) 
$$

$$
C^{2m}(2m) = L(2m) = \sum_{l=0}^{2m} (-1)^l {2m \choose l} \sum_{r=0}^{2m-l}  s(2m-l,r). \eqno(11a)
$$

{\it For the integer odds values} $ p = 5,7,9, \ldots $ {\it we have the representation }

$$
G^p(p) = L(p) = (2/e) + \sum_{k=0}^p (-1)^k {p \choose k} B(p-k). \eqno(11b)
$$ 

\section {Auxiliary results.}

{\bf 1.} {\it  In the symmetrical case for all the values}  $ p \in [4, \infty ) $ {\it we have: }

$$
K(p) = (2/e) \sum_{n=1}^{\infty} n^p \ I_n(1)  = W(p). \eqno(12)
$$

 Namely, for the values $ \tau_1, \ \tau_2 $ from (3) we receive for the values $ n=1,2,\ldots:$

$$
{\bf P}(\tau_1 - \tau_2 = n) = e^{-1} \sum_{k=0}^{\infty} \frac{2^{-k} 2^{-(n+k)}}{k! \ (k+n)!} =
I_n(1)/e.
$$

{\bf 2. } On the basis of the equality (12)  we can offer the numerical algorithm for $ K(p) $ 
investigation, calculation and estimation.  
For the improvement of speed of convergence of series (12) we can write:

$$
2 \pi \ I_n(1) = \int_{-\pi}^{\pi} \exp(\cos(\theta)) \ \cos(n \theta) \ d \theta,  
$$
 (see, for example, [16], p. 958, formula 5.)  We obtain after the integration by parts 
$$
2 \pi \ I_n(1) = (-1)^m n^{-2m} \int_{-\pi}^{\pi} (\exp(\cos \theta))^{(2 m)} \ 
\cos (n \theta) \ d \theta,
$$
$ m = 1,2, \ldots. $ Using the method of mathematical induction we conclude:

$$
(\exp \cos (\theta))^{(2m)} = \exp(\cos(\theta)) \ P_{2m}(\cos (\theta)), 
$$
where $ P_{2m}(x) $ are a polynomials of degree $ 2 m $ which may be calculated by means 
of the recursion 

$$
P_{2m+2}(x) = (1-x^2) \left(P^{//}_{2m} + 2 P^/_{2m}(x) + P_{2m}(x)  \right) -
$$
$$
x \left(P^/_{2m}(x) + P_{2m}(x) \right) 
$$
 with initial condition $ P_0(x) = 1. $  Therefore, we get the following representation for 
$ K(p): $

$$
\pi \ e \ K(p) = \sum_{n=1}^{\infty} n^{p-2m} \int_{-\pi}^{\pi} \exp \cos(\theta) \ P_{2m}(\theta) \ 
\cos (n \theta) \ d \theta.  \eqno(13)
$$
{\bf 3. Corollary. } {\it For the even numbers } $ p = 2 m, \ m = 1,2,3, \ldots $ {\it all the numbers }
$ K(p) = K(2m), L(p) = L(2m) $ {\it are integer.} \par

 In fact, it follows from formula (12) that 
$$
  K(2m) =  (\pi \ e)^{-1} \ \sum_{n=1}^{\infty} \int_{-\pi}^{\pi} g^{(2m)}(\theta) \cos(n \theta)
 \ d \ \theta = 
$$

$$
e^{-1} (\exp(\cos \theta))^{(2m)}(0) = (-1)^m P_{2m}(1).
$$

 It is easy to verify that all the coefficients of polynomials $ P_{2m}(x) $ are integer; thus, the 
number $ P_{2m}(1) $ is integer.\par
 The second conclusion of our corollary follows from the formula (10), as long as all the 
Stirling's numbers are integer. \par

{\bf 4. For example,} $  K(6) = 31, L(6) = 41. $
 For the non - integer values $ p $ we can use the method described above. We 
obtained:

\vspace{10mm}

\begin{tabular}{|c|c|c|c|c|c|}
\hline
p &  K(p) & L(p) & p & K(p) & L(p) \\
\hline
\hline
2 & 1 & 1 & 10.5 & 14000.4 & 41385.2  \\
\hline
4 & 4 & 4 & 11 & 30403.2 & 98253.7 \\
\hline
4.5 & 6.3358 & 6.6712 & 11.5 & 67091.3 & 236982 \\
\hline
5 & 10.4118 & 11.7358 & 12 & 150349 & 580317 \\
\hline
5.5 & 17.686 & 21.538 & 12.5 &  341951.2 & 1.44191E+006 \\
\hline 
6 & 31 & 41 & 13 & 788891.0 & 3.63328E+006 \\
\hline
6.5 & 55.819 & 80.5508 & 13.5 & 1.84518E+006 & 9.27951E+006 \\
\hline
7 & 103.22 & 162.7358 & 14 & 4.37346E+006 & 2.40112E+007 \\
\hline
7.5 & 192.45 & 337.176 & 14.5 & 1.04998E+007 & 6.29176E+007 \\
\hline 
8 & 379 & 715 & 15 & 2.55231E+007 & 1.66888E+008 \\
\hline
8.5 & 757.7 & 1549.28 & 15.5 & 6.27927E+007 & 4.47926E+008 \\
\hline
9 & 1126.5 & 3425.7358 & 16 & 1.56298E+008 & 1.21607 E + 009 \\
\hline
9.5 & 3015.0 & 7721.29 & 16.5 & 3.93475E+008 & 3.33839E+009 \\
\hline
10 & 6556 & 17722 & 17 & 1.00153E+009 & 9.26407E+009  \\
\hline
\end{tabular}

\vspace{10mm}

\begin{tabular}{|c|c|c|}
\hline
p & K(p) & L(p) \\
\hline
\hline
17.5 & 2.57666E+009 & 2.59791E+010 \\
\hline
18   & 6.69849E+009 & 7.36008E+010 \\
\hline
18.5 & 1.75916E+010 & 2.106E + 011 \\
\hline
19 & 4.66582E+010  & 6.08476 + 011 \\
\hline
19.5 & 1.24952E+011 & 1.77473E+012 \\
\hline
20 & 3.37789E+011 & 5.22427E+012 \\
\hline
20.5 & 9.21603E+011 & 1.55177E+013 \\
\hline
21 & 2.53714E+012  & 4.64999E+013 \\
\hline
\end{tabular}

\vspace{10mm}

{\bf 5. } Using the discrete analog of the saddle - point method ([24], p. 262 - 264), [10]),
 we find that 

$$
M(p) = [p/\log p] \cdot (1+ \varepsilon(p)), \
$$
where at $ p \to \infty $
$$
 \varepsilon(p) = \Delta(p) + \Delta^2(p) - \delta(p) \ \Delta(p) \ (1 + o(1)).  \eqno(14)
$$
Hence

$$
N(p) = [p/\log(2p)] \cdot (1+ \varepsilon(2p)) =
$$

$$
 [p/\log p] \cdot \left[1+ \Delta(p) + \Delta^2(p) - \delta(p) \ \Delta(p) \
 (1+\log 2)(1+o(1) \right].
$$

 Define for the values $ p \ge P_0 = 700 $ the following functions and constants:

$$
C_{14} = (1-\log \log P_0/\log P_0) \approx 1.402365,
$$

$$
C_{15} = 2 \cdot \left[(1+4 \Delta^2(P_0))^{1/2} + 1 \right] \approx 0.928958,
$$

$$
\zeta(p) = \log 2/\log(2p),
$$

$$
\varepsilon_+(p) = \Delta(p) + C_{14} \Delta^2(p), \ \varepsilon_-(p) = 
\Delta(p) + C_{15} \Delta^2(p),
$$

$$
M_+ = M_+(p) = [p/\log p] \cdot (1 + \varepsilon_+(p)), 
$$

$$
M_- = M_-(p) = [ p/\log p] \cdot (1+ \varepsilon_-(p)), \eqno(15a)
$$

$$
N_+(p) = [p/(e \cdot \log(2p))] \cdot(1+\varepsilon_+(2p)),
$$

$$
N_-(p) = [p/(e \cdot \log(2p))] \cdot(1+ \varepsilon_-(2p)). \eqno(15b)
$$

 More exact calculation show us that for all the values $ p \ge P_0 $
$$ 
 M_-(p) \le M(p) \le M_+(p), \ N_-(p) \le N(p) \le N_+(p).
$$

  Namely, it is very simple to see that $ \forall p \ge P_0 \ \Rightarrow $
$$
M_- \log M_- < p = M \log M < M_+ \log M_+.
$$

{\bf 6.} Let us denote 

$$
 \ b_1(x,p) =  x^p/\Gamma(x+1), 
$$
where 
$$
 \Gamma(x) = \int_0^{\infty} y^{x-1} \ e^{-y} \ dy 
$$
is the usually gamma function;

$$
  b_2(x,p) = x^p / \left(2^x \ \Gamma(x+1) \right);
$$

$$
V(x,p) =  p \log x - x \log x + x,  
$$

$$
X(p) = V(M(p),p)/p = \sup_{x \ge 4} V(x,p)/p,
$$

$$
W(x,p) = p \log x - x \log x + x(1-\log 2), 
$$

$$
Y(p) = W(N(p),p)/p = \sup_{x \ge 4} W(x,p)/p.
$$

 We have using the equality (14): $ X(p) = \log[p/(e \cdot \log p) ] + $

$$
\Delta(p) + \delta(p) + [\log(1+ \varepsilon(p)) - \varepsilon(p) + \Delta(p)
 \varepsilon(p) ] +
$$  

$$
\{\delta(p) (\varepsilon(p) - \log(1+\varepsilon(p)) \} - \delta(p) \varepsilon(p) 
\log(1+ \varepsilon(p))) = 
$$

$$
\log[p/(e \cdot \log p)] +  X_0(p), 
$$
where for $ p \ge P_0 \ X_2(p) < X_0(p) < X_1(p), \ X_1(p) \stackrel{def}{=} $ 

$$
 \Delta(p) + \delta(p) + \Delta(p) \varepsilon_+(p) + \delta(p) 
[\varepsilon_+(p) -\log(1+ \varepsilon_+(p)) ], \eqno(16a)
$$

$$
X_2(p) \stackrel{def}{=} \Delta(p) + \delta(p) + [\log(1+\varepsilon_-(p)) - \varepsilon_-(p)] -
$$

$$
-\delta(p) \varepsilon_-(p) \log(1+\varepsilon_-(p)). \eqno(16b)
$$

The function $ p \to X_1(p), \ p \in [P_0, \infty) $ is monotonically decreasing and 
$$
\exp(X_1(P_0)) < 1.7563, \ \lim_{p \to \infty} X_1(p) = 0. \eqno(16c)
$$

 At the same manner we get: $ Y(p) = \log[p/(e \cdot \log p)] + Y_0(p), $ where

$$
Y_0(p) = \log(1-\zeta(p)) -(1+\varepsilon(2p)) \times
$$
$$
 [1-\zeta(p) - \Delta(2p) + \delta(2p) \log(1+\varepsilon(2p)) ] + \delta(2p)
 (1+\varepsilon(2p)) \stackrel{def}{=}
$$

$$
\log g(p) + Y_0(p), \ Y_2(p) \le Y_0(p) \le Y_1(p), 
$$

$$
Y_1(p) \stackrel{def}{=} \Delta(2p) + \delta(2p) + (1+\varepsilon_+(2p)) \cdot \delta(p) \log 
2/(1+\delta(p) \log 2) +
$$

$$
\varepsilon_+(2p)[\Delta(2p) + \delta(2p)],  \eqno(16d)
$$

$$
 Y_2(p)\stackrel{def}{=} \Delta(2p) + \delta(2p) + \varepsilon_-(2p)
 [\Delta(2p) + \delta(2p)], \eqno(16e)
$$

where the function $ p \to Y_1(p), \ p \in [P_1, \infty), \ P_1 = 10^6 $ is monotonically 
decreasing and 

$$
\exp(Y_1(P_1)) < 1.442, \ Y_1(p) \downarrow 0, \ p \to \infty, \ \lim_{p \to \infty} Y_2(p) = 0.
\eqno(17)
$$

{\bf 7. Upper bound for } $ L(p). $ Assume in this section that $ p \ge P_0 = 700. $
 We have for the values $ p \ge P_0, $ using the well - known Stirling's formula:

$$
 e \cdot L(p) -1.5 = \sum_{n=3}^{\infty} b_1(n-1,p) \le \int_2^{\infty} b_1(x,p) \ dx + 
\sup_{x \ge 3} b_1(x,p) \le
$$

$$
(2 \pi)^{-1/2} \exp(p \cdot X(p)) + (2 \pi)^{-1/2} \int_2^{\infty} \exp (V(x,p)) \ dx.
$$
 Split the last integral into three parts  so that

$$
J(p) \stackrel{def}{=} \int_2^{\infty} \exp(V(x,p)) \ dx = J_1 + J_2 + J_3, \
 J_j = J_j(p), \ j = 1,2,3,
$$
 where 

$$
J_1(p) = \int_2^{M - \sqrt{p}} \exp(V(x,p)) \ dx, \ J_2 = \int_{M - \sqrt{p}}^{M + \sqrt{p}}
\exp (V(x,p)) \ dx, 
$$

$$
J_3 = \int_{M + \sqrt{p}}^{\infty} \exp (V(x,p)) \ dx,
$$ 
we have for the integral $ J_2, $ taking into account the inequalities $ M_- < M < M_+ $
and inequality: $ p \in [M-\sqrt{p}, M+\sqrt{p} ] \ \Rightarrow $

$$
  V(x,p) \le  p X(p) - 0.5 (x-M)^2 \cdot \left( p^2 M^{-2}_+  \right) <
$$

$$
 p X(p) - 0.5(x - M)^2 \cdot p \cdot  M_+^{-2} (p):
$$

$$
J_2 \le \exp(p \cdot X(p)) \cdot \int_{M-\sqrt{p}}^{M + \sqrt{p}} \exp 
\left(-0.5 p (x-M)^2 \ M^{-2}_+ \right) \ dx < 
$$

$$
\exp(p \cdot X(p)) \cdot \int_{-\infty}^{\infty} \exp 
\left(-0.5 p (x-M(p))^2 \  M^{-2}_+ \right) \ dx =
$$

$$
\sqrt{2 \pi} \ \exp(p \cdot X(p)) \ M_+ /\sqrt{p} \le \exp(p \cdot X(p)) \cdot \Psi_1(p), 
$$
where
$$
\Psi_1(p) = \sqrt{2 \pi \ p} \ \cdot  \left[1 + \Delta +  C_{14} \Delta^2 \right]
 / \log p.
$$

 Now we estimate  the integral $ J_3. $  For the values $ x \ge M + \sqrt{p} $ are true 
the following inequalities: 

$$ 
V(x,p) \le p X(p) -   0.5 \cdot (2 p) \cdot (p/M^2_+(p)) \le p X(p) -
$$

$$
  \log^2 p \cdot (1 + \Delta + C_{14} \Delta^2 )^{-2};
$$

$$
dV(x,p)/dx \le - p /M^2_+(p) \left[x-M(p) - \sqrt{p} \right];
$$

therefore $ J_3 \le  \exp(p \cdot X(p)) \cdot \Psi_2(p), $ where $ \Psi_2(p) = $

$$
\exp \left(-\log^2 p \cdot (1 + \Delta(p) + C_{14} \Delta^2(p))^{-2} \right) \times
$$

$$
\int_{M+ \sqrt{p}}^{\infty} \exp \left(- p \ M_+^{-2}  (x-M - \sqrt{2p})   \right) \ dx  =
$$

$$
 \exp \left(-\log^2 p \cdot (1+\Delta(p) + C_{14} \Delta^2(p))^{-2} \right) \times
$$

$$
  p \cdot \left(1 + \Delta + C_{14} \Delta^2 \right)^2 \cdot \log^{-2}(p)
$$

and analogously we find the upper estimate for $ J_1. $  \par
 Thus, $ L(p) < e^{-1} \cdot \exp(p \cdot X(p)) \times $

$$
 \left[ 1.5 \exp(-p \cdot X(p)) +
(2 \pi)^{-1/2} + \Psi_1(p) + 2 (2 \pi)^{-1/2} \ \Psi_2(p) \right] =
$$
 
$$
\exp(p \cdot X(p) ) \cdot \Psi^p_3(p), \eqno(18a)
$$
where we find by the direct calculations: $ \ \Psi_3(P_0) \le 1.00826 $ and at $ p \ge P_0 $
$$ 
\Psi_3(p) \downarrow 1, \ p \to \infty; \ \Psi_3(p) \le 1 + C_{18} \
\log p/p.  \eqno(18b)
$$

{\bf 9. Low bound for} $ L(p). $  Denote $ q = p - 1/2. $ We obtain using the Sonin's 
estimate for factorials:

$$
e L(p) \ge \sum_{n=4}^{\infty} b_1(n-1,p) = \sum_{n=3}^{\infty} b_1(n,p) \ge 
$$

$$
\int_{4}^{\infty} b_1(x,p) \ dx \ge (2 \pi)^{-1/2} \ \exp(-1/12) \ \int_{4}^{\infty} \ 
\exp(V(x,q)) \ dx \ge
$$

$$
(2 \pi)^{-1/2} \ \exp(-1/12) \ \int_{M(q)}^{M(q) + \sqrt{ q}} \exp(V(x,q)) \ dx.
$$

 Since the following implication holds: 
 $ \ q \in [M(q), M(q) + \sqrt{q} ] \ \Rightarrow $

$$
V(x,q) \ge q \ X(q)  - 0.5 (x-M(q))^2 \ q M^{-2}_-(q),
$$

we have:

$$
e L(p) \ge (2 \pi)^{-1/2} \ \exp(-1/12) \ \exp(q X(q)) \times
$$

$$
\int_{M(q)}^{M(q) + \sqrt{ q}} \exp \left[-0.5 q \ M^{-2}_-(q) \ (x-M(q))^2 \right]\ dx \ge
M_-(q) \ \times
$$

$$
0.5 \ \exp(-1/12) \ \sqrt{q} \ \exp(q X(q)) \ \left[1-\exp \left(-q^2/M^2_- \right)  \right] =
$$

$$
e \cdot \exp(p \cdot X(p)) \cdot \Psi^p_4(p),
$$
where 

$$
\Psi_4(p) \downarrow 1, p \to \infty; \ \Psi_4(p) \ge 1 + C_{19} \log p/p. 
$$
 Thus,

$$
\exp(p \cdot X(p)) \cdot \Psi^p_4(p) \le L(p) \le \exp(p \cdot X(p)) \cdot \Psi^p_3(p), \eqno(19a)
$$

$$
\Psi_3(p) \le 1 + C_{19} \log p/p, \ \Psi_4(p) \ge 1 + C_{20} \log p/p, \eqno(19b)
$$

$$
\Psi_{3}(p) \downarrow 1, p \to \infty; \Psi_3(P_0) \le 1.00826. \eqno(19c)
$$

{\bf 10. Upper and low bounds for } $ K(p)  $ are provided analogously to the upper 
bound for $ L(p), $ but we assume in this section that $ p \ge P_1 = 10^6. $ Briefly:

$$
\sum_{k=0}^{\infty} \frac{4^{-k}}{k! \ (n+k)!} < \frac{1}{n!} \sum_{k=0}^{\infty} 
\frac{4^{-k}}{k!} = \frac{\sqrt[4]{e}}{n!} < \frac{ 1.285}{n!},
$$
 hence 
$$
K(p) < 2 e^{-3/4} \sum_{n=1}^{\infty} \frac{n^p \ 2^{-n}}{n!} = 2 \sqrt[4]{e} \cdot B(0,p,1/2).
$$

 Further, we conclude, again using the  Stirling estimate for factorials:
$$ 
0.5 \ e^{3/4} \ K(p) = \sum_{n=1}^{\infty} b_2(n,p) \le \int_2^{\infty} b_2(x,p) \ dx +
\sup_{x \ge 2} b_2(x,p) \le
$$

$$
 (2 \pi)^{-1/2} \exp(p \cdot Y(p)) + (2 \pi)^{-1/2} \int_{2}^{\infty} \exp (W(x,p)) \ dx. 
$$
 We have, again split the last integral:

$$
 I_4 \stackrel{def}{=}  \int_2^{\infty}\exp(W(x,p)) \ dx =
$$
$$
 \left( \int_2^{N(p) - \sqrt{p}}
+ \int_{N(p)-\sqrt{p}}^{N(p) + \sqrt{p}} + \int_{N(p) + \sqrt{p}}^{\infty} \right) \exp W(x,p) dx =
$$

 $ I_5 + I_6 + I_7. $ As long as at $ x > N(p) + \sqrt{p} \ \Rightarrow $

$$
W(x,p) \le p Y(p) - 0.5p^2 N^2_+(p) = p Y(p) - 0.5 \log^2 p \cdot (1+\varepsilon_+(2p))^{-2},
$$

$$
dW/dx \le -p N^{-2}_+(x-N-\sqrt{p}), 
$$
we obtain:

$$
I_7 \le \exp( p Y(p)) \cdot  p \ \log^{-2} p \ (1+\varepsilon_+(2p))^2 \times
$$

$$
 \exp \left(-0.5 \log^2 p \ (1+\varepsilon_+(2p))^{-2} \right).
$$

 Further, if $ x \in [N(p) - \sqrt{p}, N(p) + \sqrt{p}], $ then 

$$
W(x,p) \le p Y(p) -0.5 p N^{-2}_+(p) \cdot(x-N(p))^2.
$$
 Therefore
$$
I_6 \le \exp(p Y(p)) \cdot \sqrt{p} \ (1+\varepsilon_+(2p))/\log (2p)
$$
and $ K(p) \le 2 e^{-3/4} \ \exp(p Y(p)) \times $

$$
 \left[(2 \pi)^{-1/2}  + \sqrt{p} (1+\varepsilon_+(2p)) /\log(2p) + 2 (2 \pi)^{-1/2} p 
\ \log^{-2}p \right] \times
$$

$$
\left[(1+\varepsilon_+(2p))^2 \cdot \exp(-0.5 \log^2 p \ (1+\varepsilon_+(2p))^{-2} \right].
\eqno(20)
$$

 {\it Low bound for } $ K(p). $ We have:  $ 0.5 \ e \ K(p) > $ 

$$
 \sum_{n=1}^{\infty} n^p 2^{-n} /n! = B(0,p,1/2) > \exp(-1/12) \ (2 \pi)^{-1/2} 
\ \exp(q Y(q)) \times
$$

$$
\int_{N(q)}^{N(q) + \sqrt{q} } \exp \left[ -0.5 \ q \ N^{-2}_-(q) \ (x-N(q))^2  \right] \ dx \ge
$$

$$
\exp(-1/12) \ \sqrt{\pi/2} \ q^{-1/2} \ N_-(q) \ \left(1-\exp \left(-q^2/N^2_-(q) \right) \right)
/\log (2q) =
$$

$$
\exp(-1/12) \ \sqrt{\pi/2} \ \exp(q Y(q)) \ \sqrt{q} \ (1+ \varepsilon_-(2q)) \times
$$

$$
\left( 1-\exp \left(-q^2/N^2_-(q) \right) \right)/\log(2q).
$$

 Further estimations are like to the estimation for $ L(p) $ and may be omitted. Result:

$$
\exp(p \cdot Y(p)) \cdot \Psi_6^p(p) \le K(p) \le \exp(p \cdot Y(p)) \cdot \Psi^p_5(p), 
\eqno(21a)
$$
where at $ p \ge P_1 $ 

$$
\Psi_{5}(p) \le 1 + C_{19} \log p /p, \ \Psi_6(p) \ge 1 + C_{20} \log p/p; \eqno(21b)
$$

$$
\Psi_{5}(p) \downarrow 1, \ p \to \infty;  \ \Psi_5(P_1) \le 1.000833. \eqno(21c)
$$

{\bf 11}. For the correct calculations (by computer) we need to estimate the derivatives 
of our functions $ L(p), K(p). $  We show here the estimation of  derivatives
$ L^{(m)}(p), m = 1,2 \ldots. $ Namely,  $ \ e \cdot L^{(m)}(p) = $ 
$$
\sum_{n=3}^{\infty} \frac{(n-1)^p \ \log^m(n-1)}{n!}  \le \sum_{n=3}^{\infty} 
\frac{(n-1)^p}{(n-1)!} \cdot \frac{\log^m n}{n} <
$$

$$
\sum_{n=2}^{\infty} \frac{n^p}{n!} \cdot \left( \sup_{n \ge 3} \frac{\log^m n}{n}  \right) =
$$

$$
\left(\frac{m}{e} \right)^m \cdot \sum_{n=2}^{\infty} \frac{n^p}{n!} = \left(\frac{m}{e}\right)^m 
\cdot (e B(p) -1). \eqno(22)
$$
 The derivative $ K^{(m)}(p), \ m = 1,2,\ldots $ we estimate analogously.\par
 It follows from this estimation that the functions $ L(p) $ and $ K(p) $ are infinitely 
differentiable at the interval $ p \in (4, \infty). $ As long as $ L(4-0) = L(4+0) = 4, \ 
K(4-0) = K(4+0) = 4, $ both the functions $ K(\cdot), \ L(\cdot) $ are continuos in the 
semiclosed interval $ [2, \infty). $ But 

$$
\frac{ dK}{dp}(4-0) = \frac{dL}{dp}(4-0) \approx 3.149195, 
\ \frac{dK}{dp} (4+0) \approx 3.51934,  
$$

$$
\frac{dL}{dp}(4+0) \approx 3.86841, 
$$ 
therefore both the functions $ K(\cdot), \ L(\cdot) $ are not 
continuous differentiable in the set $ (2, \infty). $ 
 In the open intervals $ (2,4) $ and $ (4,\infty) $ all the functions $ L(p), K(p), C(p), S(p) $ 
are infinitely differentiable (see (22) and [10], [17] ). \par

\section{ Proof of the main results.}

{\bf Proof of theorem 1.} We find by the direct calculations that \\ $ G(C_4)/g(C_4) \approx 
1.77638, $ but for we conclude  from (17) that for the values $ p \ge P_0 = 700 $

$$
G(p) /g(p) \le 1.00826 \cdot 1.75913 = 1.77366, 
$$
hence 

$$
\argm_{p \in [4, \infty)} G(p)/g(p) \in [4, 700].
$$
 We obtain by direct calculations using usually numerical methods and by means of computer:
$$
\max_{p \in [4,700]} G(p)/g(p) = G(C_4)/g(C_4) \approx 1.77638.
$$
 Further, 

$$
\inf_{p \ge 4 } G(p)/g(p) = \min \left\{\min_{p \in [4, 700]} G(p)/g(p), \ \inf_{p > 700}
G(p) /g(p) \right\}.
$$
We obtain by computer calculations that 

$$
\min_{p \in [4,700]} G(p)/g(p) \approx 1.332,
$$
and it follows from (18a), (18b), (18c) and (19) that 

$$
\inf_{p > 700} G(p)/g(p) = \lim_{p \to \infty} G(p)/g(p) = 1.
$$
 Thus, 
$$
\inf_{p \ge 4} G(p)/g(p) = \lim_{p \to \infty} G(p)/g(p) = 1.
$$
 Analogously, $  S(C_{10})/g(C_{10}) \approx 1.53572, $ but we have for the values 
$ p \ge P_1 $

$$
S(p)/g(p) \le 1.0008333 \cdot 1.443 < 1.4444.
$$
 Therefore 
$$
\argm_{p \in [4, \infty)} S(p)/g(p) \in [4, 1000000].
$$
We obtained after some technical difficulties using computer:

$$
\max_{p \in [4, P_1]} S(p)/g(p) = S(C_{10})/g(C_{10}) \approx 1.53572.
$$
 Further, 

$$
\inf_{p \ge 4} S(p)/g(p) = \min \{ \min_{P \in [4, P_1]} S(p)/g(p), \inf_{p \ge P_1}
S(p)/g(p) \} =
$$

$$
\lim_{p \to \infty} S(p)/g(p) = 1.
$$
 The another assertions of theorem 1 are obtained analogously. \par
{\bf Proof of theorems 2 and 3.} It follows from inequalities (19a), (19b), (19c) and (21a), 
(21b), (21c) that

$$
\exp(X(p)) \cdot (1+C_{20} \log p/p) \le G(p) \le \exp(X(p)) \cdot (1+C_{19} 
\log p/p), \eqno(23a)
$$

$$
\exp(Y(p)) \cdot (1+C_{20} \log p/p) \le S(p) \le 
\exp(Y(p)) \cdot (1+ C_{19} \log p/p). \eqno(23b)
$$
 Substituting the expression (16a) and (16b) into the formula (23a), we obtain after 
some simple calculation  our assertions (9b); (9a) provided analogously. \par
 Finally, substituting expressions (16d,e)  into (23b), we obtain (7a), (7b). \par

{\bf Proof of theorem 4.}  Since $ {\bf E} \theta_{(r)} = 1, r=0,1,2, \ldots, $ 
(see (4)), we conclude

$$
{\bf E} \theta^k = {\bf E} \sum_{l=0}^k s(k,l) \theta_{(r)} =  \sum_{l=0}^k  s(k,l).
$$
 Formula (11a) follows from the  binomial formula. Equality (10)  proved analogously. \par
Let us prove (11b). Since 

$$
\sum_{n=2}^{\infty} n^p / n! = e \cdot B(p) - 1, \ p = 1,2,3, \ldots; \
\sum_{n=2}^{\infty} 1/n! = e - 2, 
$$
we conclude for the values $ p = 5,7,9, \ldots:  \ e \cdot L(p) = $

$$
 1 + \sum_{n=2}^{\infty} \frac{(n-1)^p}{n!} =  1+ \sum_{n=2}^{\infty} (n!)^{-1} \cdot
\sum_{k=0}^p (-1)^k {p \choose k}  n^{p-k}  = 
$$

$$
1+ \left[ \sum_{k=1}^p (-1)^k {p \choose k}
(e B(p-k) - 1) \right] -
$$

$$
(e B(0) - 2) = 2 + e \cdot \sum_{k=0}^p (-1)^k {p \choose k} B(p-k) - 
\sum_{k=0}^p (-1)^k {p \choose k}  = 
$$

$$
2+ e \sum_{k=0}^p (-1)^k {p \choose k}  B(p-k).
$$

\section{Concluding Remark.} 

 Our results allow us to obtain some generalizations on the Hilbert space symmetrical 
distributed random variables $ \{ \eta(i) \}, \ i = 1,2,3, \ldots. $ Let $ \ (H, |||\cdot|||) $  
be a separable Hilbert space with the norm $ |||\cdot|||, \ {\bf P}(\eta(i) \in H) = 1,
\ \forall i = 1,2,3 \ldots \ ||\eta||_p \stackrel{def}{=} {\bf E}^{1/p} \left(|||\eta(i)|||^p \right) 
< \infty, \ p \ge 4, $ 

$$
Z(p) = \sup_{ \{\eta(i)\}} \sup_n \frac{||\sum \eta(i)||_p}
{\max \left( ||\sum \eta(i)||_2, (\sum ||\eta(i)||_p^p)^{1/p} \right)}.
$$
 Utev ([17], [18]) proved that  $ Z(p) = S(p), \ p \ge 4 $ (in our notations). Therefore 
$$
 1 = \inf_{p \ge 4} Z(p)/g(p) < \sup_{p \ge 4} Z(p)/g(p) = C_9 \approx 1.53572,
$$ 

$$
1 = \inf_{p \ge 15} Z(p)/h(p) < \sup_{p \ge 15} Z(p)/h(p)  = C_{11} \approx 1.03734.
$$

 Probably, it is interest to obtain the exact constants in the moment inequalities 
for sums of independent nonnegative random variables in the spirit of articles 
[7], [8], [12] etc.\\

\vspace{2mm}

{\bf Aknowledgements}. I am very grateful to prof. V.Fonf, M.Lin (Ben Gurion University, 
Beer - Sheva, Israel) for useful support of these investigations. \par
 I am grateful to prof. G.Schechtman (Weizman Institute of Science, Rehovot, Israel) for 
attention.\par
 This investigation was partially supported by ISF(Israel Science Foundation), grant 
$ N^o $ 139/03.\par

\newpage

{\bf References.} \\

1. William Feller. An Introduction to Probability Theory and its Applications. Volume 2,
1966. John Wiley and Sons Inc.,  New York.\\
2. Petrov V.V.  Limit theorems of Probability Theory.  Sequences of Independent  Random 
Variables. 1995. Oxford Science Publications. Claredon Press. Oxford, UK. \\

3. Dmitrovsky V.A., Ermakov V.A, Ostrovsky E.I. The Central Limit Theorem for weakly dependent 
Banach - Valued Variables. Theory Probab. Appl., 1983, V. 28, 1, p. 89 - 104.\\

4. Talagrand M. Isoperimetry and integrability of the sum of independent Banach - space valued 
random Variables. Ann. Probab., 1989, V. 17 p. 1546 - 1570.\\

5. Buldygin V.V., Mushtary D.M., Ostrovsky E.I., Puchalsky A.I. New Trends in Probability Theory
and Statistics. 1992, Springer Verlag, New York - Berlin - Heidelberg - Amsterdam.\\

6. Dharmadhikari S., Jogdeo K. Bounds on the moment of certain random variables. Ann. Math. 
Statist., 1969, V. 40, 4, p. 1506 - 1518.\\

7. Rosenthal H.P. On the subspaces of $ L^p \ (p > 2) $ spanned by sequences of independent 
variables. Israel J. Math., 1970, $ N^o 3, $ p. 273 - 303.\\

8. Johnson W.B., Schechtman G., Zinn J. Best constants in moment inequalities for linear 
combination of independent  and changeable random variables. Ann. Probab., 1985, V. 13 p. 234 - 253.\\

9. Figel T., Hitczenko P., Johnson W.B., Schechtman G., Zinn J. Extremal properties of Rademacher 
functions with applications to the Khintchine and Rosenthal Inequalities. Trans. of the  Amer.
Math. Soc., 1997. V. 349, p. 997 - 1027.\\

10. Ibragimov R., Sharachmedov Sh. On the exact constant in the Rosenthal Inequality.
 Theory Probab. Appl., 1997, V. 42.  p. 294 - 302.\\

11. Ibragimov R., Sharachmedov Sh. The exact constant in the Rosenthal Inequality for sums random
 variables with mean zero. Probab. Theory Appl., 2001, V. 46, 1, p. 127 - 132.\\

12.Latala R. Estimation of moment of sums of independent real Random variables. Ann. Probab.,
 1997, V. 25, 3, p. 1502 - 1513.  \\

13. Peshkir G., Shiryaev A.N. The Khintchin's inequalities and martingale expanding  sphere of 
their Actions. Russian . Math. Surveys, 1995, V. 50. 2(305), p. 849 - 904.\\

14. Kennet H. Rosen (Editor - in - Chief), John G. Michaels, Johnatan L. Gross, Jerrold W. Grossman,
 Douglas R. Shier. Handbook of Discrete and Combinatorial Mathematics. 2000. CRC Press, Boca Raton.
London, New York, Washington. \\

15. Fedorjuk M.V. The Saddle - point method. Kluwner, 1990. Amsterdam, New York.\\

16. Gradstein I.S., Ryszik I.M. The Tables of Integrals, Sums and Products. 1980, Academic Press,
New York, London, Toronto, Sydney, San Francisco.\\

17. Utev S.A. The extremal problems in the moment Inequalities. In: Coll. Works of Siberian 
Branch of Akademy Science USSR, Limit theorems in Probability Theory; 1985, V. 23, p. 56 - 75.\\ 

18. Utev S.A. The extremal Problems in Probability Theory. Theory Probab. Appl., 1984, V. 28, 2,
p. 421 - 422.\\

19. Ostrovsky E.I. Exponential estimations for Random Fields and their Applications (in Russian).
 Obninsk, OINPE, 1999.\\

20. Pinelis I.F and Utev S.A. Estimates of moment of sums of independent random variables.
Theory Probab. Appl., 1984, V. 29, p. 574 - 578.\\

21. Satchkov V.N. Combinatorial Methods in Discrete Mathematics. Cambridge University Press,
 1996, Cambridge, UK. \\

22. Johnson W.B., Maurey B., Schechtman G., Tzafriri L. Symmetric structures in Banach Spaces.
 Memoirs of the A.M.S. 1979, {\bf 217.} \\












\end{document}